

\baselineskip=14pt
\parskip=10pt
\def\halmos{\hbox{\vrule height0.15cm width0.01cm\vbox{\hrule height
  0.01cm width0.2cm \vskip0.15cm \hrule height 0.01cm width0.2cm}\vrule
  height0.15cm width 0.01cm}}
\font\eightrm=cmr8 
\font\eighttt=cmtt8
\magnification=\magstephalf

\def\1{{\overline{1}}}
\def\2{{\overline{2}}}
\parindent=0pt
\overfullrule=0in

\def\frac#1#2{{#1 \over #2}}
\bf
\centerline
{
Alexander Burstein's Lovely Combinatorial Proof of John Noonan's Beautiful Theorem
}
\centerline
{
that the number of $n$-permutations that contain the Pattern 321 Exactly Once Equals}
\centerline
{
 (3/n)(2n)!/((n-3)!(n+3)!)
}
\rm
\bigskip
\centerline{ {\it
Doron 
ZEILBERGER}\footnote{$^1$}
{\eightrm  \raggedright
Department of Mathematics, Rutgers University (New Brunswick),
Hill Center-Busch Campus, 110 Frelinghuysen Rd., Piscataway,
NJ 08854-8019, USA.
{\eighttt zeilberg  at math dot rutgers dot edu} ,
\hfill \break
{\eighttt http://www.math.rutgers.edu/\~{}zeilberg/} .
Based on a mathematics colloquium talk at Howard University, Oct. 14, 2011, 4:10-5:00pm,
where Alex Burstein was presented with  a \$25 check prize promised in John Noonan's
1996 Discrete Math paper. Oct 18., 2011.
Exclusively published in the Personal Journal of Shalosh B. Ekhad and Doron Zeilberger 
{\eighttt http://www.math.rutgers.edu/\~{}zeilberg/pj.html} and {\eighttt arxiv.org} .
Supported in part by the NSF.
}
}

Alex Burstein[1] gave a lovely combinatorial
proof of John Noonan's[2]
lovely  theorem that the number of $n$-permutations that contain the pattern 321 exactly once
equals $\frac{3}{n}{{2n} \choose {n+3}}$.
Burstein's proof can be made even shorter as follows. Let $C_n:=(2n)!/(n!(n+1)!)$ be the Catalan numbers.
It is well-known (and easy to see) that $C_{n}=\sum_{i=0}^{n-1} C_iC_{n-1-i}$.
It is also well-known (and fairly easy to see) that the number of $321$-avoiding $n$-permutations
equals $C_n$.

Any $n$-permutation, $\pi$, with exactly one $321$ pattern can be written as $\pi_1 c \pi_2 b \pi_3 a \pi_4$,
where $cba$ is the unique $321$ pattern (so, of course $a<b<c$). All the entries to the left
of $b$, except $c$, must be smaller than $b$, and all the entries to the right of $b$, except for $a$,
must be larger than $b$, or else another $321$ pattern would emerge. 
Hence $\sigma_1:=\pi_1 b \pi_2 a$ is a $321$-avoiding permutation of $\{1, \dots , b\}$ that does not end with $b$ and
$\sigma_2:=c \pi_3 b \pi_4$ is  a $321$-avoiding permutation of $\{b, \dots , n \}$ that does not start with $b$.
This is a bijection between the Noonan set and the set of pairs $(\sigma_1, \sigma_2)$ as above (for some $2 \leq b \leq n-1)$.
For any $b$ the number of possible $\sigma_1$ is $C_b-C_{b-1}$.
Similarly, the number of possible $\sigma_2$ is $C_{n-b+1}-C_{n-b}$. Hence
the desired number is 
$$
\sum_{b=2}^{n-1} (C_b-C_{b-1})(C_{n-b+1}-C_{n-b})=
\sum_{b=1}^{n} (C_b-C_{b-1})(C_{n-b+1}-C_{n-b})
$$
$$
=\sum_{b=1}^{n} C_bC_{n-b+1}
-\sum_{b=1}^{n} C_{b}C_{n-b}
-\sum_{b=1}^{n} C_{b-1}C_{n-b+1}
+\sum_{b=1}^{n} C_{b-1}C_{n-b}
$$
$$
=C_{n+2}-2C_{n+1} - 2(C_{n+1}-C_n)+C_n=
C_{n+2}-4C_{n+1}+3C_n=\frac{3}{n}{{2n} \choose {n+3}} \quad . \quad \halmos 
$$

{\bf References}

1. Alex Burstein, {\it A short proof for the number of permutations containing the pattern 321 exactly
once},  Elec. J. Comb. {\bf 18(2)}(2011), \#P21.

2. John Noonan, {\it The number of permutations containing exactly one increasing subsequence of length three},
Discrete Math. {\bf 152}(1996), 307-313 .

\end